\documentclass[11pt,leqno]{article}
\usepackage{amsmath,amssymb,amsfonts}
\usepackage{graphicx}

\newtheorem{theorem}{{\sc Theorem}}
\newcommand{\bt}{\begin{theorem}}
\newcommand{\et}{\end{theorem}}
\newcommand{\newsection}[1]{\setcounter{equation}{0} \setcounter{theorem}{0}
\section{#1}}

\newcommand{\NI}{\noindent}
\newcommand{\bea}{\begin{eqnarray}}
\newcommand{\eea}{\end{eqnarray}}

\def \spec#1 {\mathop{#1}}

\def \b #1 {\bf #1}

\newcommand {\CC}{\centerline}

\newcommand{\ity}{\infty}
\newcommand{\raro}{\rightarrow}

\newcommand{\vsp}{\vskip 1em}

\newcommand{\be}{\begin{equation}}
\newcommand{\ee}{\end{equation}}
\newcommand{\ben}{\begin{eqnarray*}}
\newcommand{\een}{\end{eqnarray*}}

\begin{document}
\CC{\bf{Nonparametric Estimation of Linear Multiplier for SDEs} }
\CC{\bf{Driven by Generalized Stochastic Volatility}}
\vsp
\CC{B.L.S. Prakasa Rao}
\CC{CR RAO Advanced  Institute of Mathematics, Statistics}
\CC{and Computer Science, Hyderabad, India}
\CC{(e-mail address: blsprao@gmail.com)}
\vsp
\NI{\bf Abstract:} We study the problem of nonparametric estimation of the linear multiplier function $\theta(t)$ for processes satisfying stochastic differential equations of the type
$$dX_t= \theta(t)X_t\; dt+ \epsilon\; \sigma(t,X_t,Y_t) \;dW_t, X_0=x_0, 0 \leq t \leq T,$$
based on the observation of the process $\{X_t,0\leq t \leq T\},$ and investigate the asymptotic properties of the estimator as the small noise $\epsilon$ tends to zero. Here $\{W_t, t\geq 0\}$ is the standard Brownian motion and $\{Y_t, t\geq 0\}$ is a process adapted to the filtration generated by the Brownian motion. 
\vsp
\NI{\bf Keywords :} Nonparametric estimation; Linear multiplier; Kernel method; Generalized stochastic volatility; Brownian motion.
\vsp
\NI{\bf Mathematics Subject Classification :} Primary 62G05.
\newsection {Introduction} 
We now discuss the problem of estimating the function $\theta(t), 0 \leq t \leq T$ (linear multiplier) based on the observations of the  process $\{X_t, 0\leq t \leq T\}$ satisfying the stochastic differential equation
\be
dX_t=\theta(t)X_t\; dt + \epsilon \; \sigma(t,X_t,Y_t)\;dW_t, X_0=x_0, 0 \leq t \leq T  
\ee
where $\{W_t, t \geq 0\}$ is the standard Brownian motion. As pointed out by Khlifa et al. (2016,2018), such stochastic differential equations (SDEs) are used as models of a financial market in stochastic or mathematical finance. Hull  and White (1987) applied such a model for the pricing of options on assets with stochastic volatilities where $\sigma (t,x,y)=xy$ and the process $Y$ is the square root of the geometric Brownian motion. A similar model was proposed by Heston (1993) with applications to model bond and currency options with the volatility governed by the Ornstein-Uhlenbeck process. Fouque et al. (2000a,b) used the model where the stochastic volatility is driven by the Cox-Ingersoll-Ross process. Kuchuk-Iatsenko and Mishura (2015a,b) considered the case where $\sigma(t,x,y)=x \sigma_1(y)$ and $Y$ is the Ornstein-Uhlenbeck process. Parameter estimation in stochastic volatility models is studied in Ait-Sahalia and Kimmel (2007). We investigate the nonparametric estimation of the function $\theta(t)$ when the generalized stochastic volatility is present. We will denote any positive constant by C through out the sequel and the constant $C$ may differ from one line to another.
\vsp
Khlifa et al. (2016) considered SDEs of the form
$$dX_t= \theta a(t,X_t) dt + \sigma_(t,X_t)\sigma_2(t,Y_t)dW_t$$
with multiplicative stochastic volatility and considered the problem of estimation of the parameter $\theta$ through the method of maximum likelihood and obtained sufficient conditions for the strong consistency of this estimator. A similar investigation has been carried out by Khlifa et al. (2018) for SDEs of the type
$$
dX_t=\theta a(t,X_t) dt + \sigma(t,X_t,Y_t)dW_t, X_0=x_0, 0 \leq t \leq T. 
$$
Khlifa et al. (2018) provide sufficient conditions for the existence and uniqueness for the strong solution of the equation given above. Detailed result is stated in the appendix. Our aim in this paper is to estimate the drift function $\theta(t) X_t$ and the function $\theta(t)$ in the SDE (1.1), based on the observation of the process $\{X_t,0\leq t \leq T\},$ and study their asymptotic properties as the noise $\epsilon$ tends to zero. A special case of this problem was studied in Kutoyants (1994) when the function $\sigma(t,x,y)$ is identically equal to one. Our approach to the problem is similar to that of Kutoyants (1994).
\vsp
We assume that the process $Y$ is a continuous process adapted to the filtration generated by the Brownian motion $W.$ In addition, we assume that following conditions  hold for the functions $\theta(.), \sigma(t,x,y)$ which imply that the equation (1.1)  has a unique strong solution $\{X_t,0\leq t \leq T\}$ by the theorem due to Khlifa et al. (2018) stated in the Appendix.
\vsp
\NI{$(A_1)$ The trend coefficient $\theta(t),$ over the interval $[0,T],$ is bounded by a constant $L$;\\
\NI{$(C1^\prime)$ there exists a constant $K>0$ such that
$$|\sigma(t,x,y)|^2\leq K^2(1+|x|^2), t \geq 0, x,y \in R;$$
\NI{$(C2^\prime)$ there exist  constants $K_N>0, C_N>0$ and a positive function $\phi(t,y)$such that
$$|\sigma(t,x_1,y)-\sigma(t,x_2,y)|\leq K_N|x_1-x_2| \phi(t,y), t\geq 0, x_1,x_2,y \in R$$
where
$$\sup_{s\geq 0}\sup_{|x|\leq N}|\phi(s,x)|\leq C_N.$$
\vsp
Consider the differential equation in the limiting system of (1.1) as $\epsilon \raro 0,$ that is, for $\epsilon=0,$ given by
\be
\frac{dx_t}{dt}=\theta(t) x_t , x_0, 0 \leq t \leq T.
\ee
Observe that
$$x_t=x_0 \exp \{ \int^t_0 \theta(s) ds\}, 0\leq t \leq T.$$
\vsp
\noindent{\bf Lemma 1.1.} Suppose the condition $(A_1), (C1^\prime)$ and $(C2^\prime)$ hold and  let $\{X_t, 0\leq t \leq T\}$ and $\{x_t, 0\leq t \leq T\}$ be the solutions of the equations (1.1) and (1.2) respectively. Then, with probability one,
\be
(a) \;|X_t-x_t| < e^{L t} \epsilon \sup_{0\leq u\leq t}|\int_0^u \sigma(s,X_s,Y_s)dW_s|,0\leq t \leq T
\ee
and there exists a constant $C>0$ such that 
\be
(b)\;\sup_{0 \leq t \leq T} E(X_t-x_t)^2  \leq C e^{2L T} T \epsilon^2 .
\ee
\vsp
\NI{\bf Proof of (a) :} Let $a(t,X_t)= \theta(t) X_t$ and $u_t=|X_t-x_t| $. Then, by the condition $(A_1)$, it follows that
\bea
u_t & \leq & \int^t_0 \left| a(v,X_v)- a(v,x_v) \right| dv + \epsilon \;|\int_0^t \sigma(s,X_s,Y_s)dW_s|\\\nonumber
& \leq & L \int^t_0 u_v dv + \epsilon \sup_{0\leq u \leq t}|\int_0^u \sigma(s,X_s,Y_s)dW_s|.
\eea
Applying the Gronwall lemma (cf. Lemma 1.11, Kutoyants (1994), p.25),  it follows that
\be
u_t \leq  \epsilon \sup_{0\leq u \leq t}|\int_0^u \sigma(s,X_s,Y_s)dW_s| e^{Lt}. \\
\ee
\vsp
\NI{\bf Proof of (b) :} Let
$$V(t)=\int_0^t \sigma(s,X_s,Y_s)dW_s, 0\leq t \leq T.$$
Observe that 
$$E[V^2(t)]= \int_0^t E[\sigma(s,X_s,Y_s)]^2\;ds $$
and the last term is finite by the conditions $(A_1),(C1^\prime)$ and the remarks stated in the Appendix. The process $\{V(t),0\leq t \leq T\}$ is a martingale adapted to the filtration generated by the Brownian motion $W$  and it follows that
$$E[\sup_{0\leq t \leq T}V^2(t)]\leq 4 E[V^2(T)]$$ 
(cf. Klebaner (2012), Theorem 7.31, p.203). 
\vsp
From the equation (1.3), we have ,
\bea
E(X_t-x_t)^2 & \leq & e^{2Lt} \epsilon^2 E (\sup_{0\leq v\leq t}\int_0^v \sigma(s,X_s,Y_s) dW_s)^2 \\\nonumber
&\leq & 4 e^{2Lt} \epsilon^2\int_0^tE(\sigma(s,X_s,Y_s))^2 ds\\\nonumber
&=& e^{2Lt} C \epsilon^2\int_0^t E(1+ |X_s|^2)ds\;\;\mbox{(by the conditions}(A_1), (C1)) \\\nonumber
&=& e^{2Lt} C t \epsilon^2\sup_{0\leq s \leq T}E(1+ |X_s|^2)\\\nonumber
& = & e^{2Lt}\epsilon^2 C t E(1+X_0^2) \;\; \mbox{(by the remarks made in the Appendix)}\\\nonumber
\eea
for some positive constant C. Hence
\be
\sup_{0 \leq t \leq T} E (X_t-x_t)^2  \leq C \; e^{2LT} T\epsilon^2 . \\
\ee
\vsp
\newsection{Main Results}
Let $\Theta_0(L)$ denote the class of all functions $\theta(.)$ with the same bound $L$. Let
$\Theta_k(L) $ denote the class of all functions $\theta(.)$ which are uniformly bounded by the same constant $L$
and which are $k$-times differentiable with respect to $t$ satisfying the condition
$$|\theta^{(k)}(x)-\theta^{(k)}(y)|\leq L_1|x-y|, x,y \in R$$
for some constant $L_1 >0.$ Here $g^{(k)}(x)$ denotes the $k$-th derivative of $g(.)$ at $x$ for $k \geq 0.$ If $k=0,$ we interpret the function $g^{(0)}(x)$ as $g(x).$
\vsp
Let $G(u)$ \ be a  bounded function  with compact support $[A,B]$ with $A<0<B$ satisfying the condition
\vsp
\NI{$(A_2)$} $\int^B_A G(u) du =1.$
\vsp
It is obvious that the following properties are satisfied by the function $G(.):$

$$\int^\infty_{-\infty} |G(u)|^2 du < \infty;$$ 
and
$$\int^\infty_{-\infty} |u^{k+1} G(u)|^2 du <\infty, k \geq 0.$$

We define a kernel type estimator $\hat \theta_t$ of the function $\theta(t)$ by the relation
\be
\widehat{\theta}_t X_t= \frac{1}{\varphi_\epsilon}\int^T_0 G \left(\frac{\tau-t}{\varphi_\epsilon} \right) d X_\tau
\ee
where the normalizing function  $ \varphi_\epsilon \rightarrow 0 $ as \ $ \epsilon \rightarrow 0. $ Let $E_\theta(.)$ denote the expectation when the function $\theta(.)$ is the linear multiplier.
\vsp
\NI{\bf Theorem 2.1:}  {\it Suppose that the linear multiplier $\theta(.) \in \Theta_0(L)$ and  the function \ $ \varphi_\epsilon \rightarrow 0$  and $\epsilon\varphi_\epsilon^{-1}\raro 0$ as $\epsilon \raro 0.$ Suppose the conditions $(A_1),(A_2),(C1^\prime)$ and $(C2^\prime)$ hold.Then, for any $ 0 < a \leq b < T,$ the estimator $\hat \theta_t X_t$ is uniformly consistent, that is,}
\be
\lim_{\epsilon \rightarrow 0} \sup_{\theta(.) \in \Theta_0(L)} \sup_{a\leq t \leq b } E_\theta ( |\hat \theta_t X_t-\theta(t) x_t|^2)= 0.
\ee
\vsp
In addition to the conditions $(A_1),(A_2)$, suppose the following condition holds:\\
\NI{$(A_3)$} $\int_{-\infty}^{\infty}u^j G(u)du=0, j=1,2,\dots,k.$

\NI{\bf Theorem 2.2:} {\it Suppose that the function $ \theta(.) \in \Theta_{k+1}(L)$ for some $k\geq 1$ and  the conditions $(A_1)-(A_3), (C1^\prime)$ and $(C2^\prime)$ hold. Further suppose that $ \varphi_\epsilon = \epsilon^{\frac{1}{k+2}}.$ Then,}
\be
\limsup_{\epsilon \rightarrow 0} \sup_{\theta(.) \in \Theta_{k+1}(L)}\sup_{a \leq t \leq b} E_\theta (| \hat\theta_t X_t - \theta(t)x_t|^2)
\epsilon^{-\frac{2(k+1)}{k+2}} \ < \infty.
\ee
\vsp
\NI{\bf Theorem 2.3:} {\it Suppose that the function $\theta(.) \in \Theta_{k+1}(L)$ for some $k\geq 1$ and  the conditions $(A_1)-(A_3), (C1^\prime)$ and $(C2^\prime)$ hold. Further suppose that 
$\varphi_\epsilon= \epsilon^{\frac{1}{k+2}}.$ Let $J(t)= \theta(t)x_t.$ Then, as $\epsilon \raro 0,$  the asymptotic distribution of the random variable
$$ \epsilon^{-\frac{2k+3}{2k+4}} (\hat\theta_t X_t-J(t)- \frac{J^{(k+1)}(t)}{(k+1) !} \int^\infty_{-\infty} G(u) u^{k+1}\ du)$$
is  the limiting distribution of the random variable
$$\frac{1}{\varphi_\epsilon}\int_0^T G \left(\frac{\tau-t}{\varphi_\epsilon} \right)\sigma(\tau,X_\tau,Y_\tau)dW_\tau$$
as $\epsilon \raro 0$ with mean zero and the variance
}
 $$\nu(t)\int_R G^2(u)du$$
where
$$ \nu(t)=E[\sigma(t,X_t,Y_t)]^2.$$ 
\vsp
\NI{Remarks 2.1:} The limiting distribution in Theorem 2.3 is Gaussian if 
$$\frac{1}{\varphi_\epsilon^2}\int_0^T G^2 \left(\frac{\tau-t}{\varphi_\epsilon} \right)\sigma^2(\tau,X_\tau,Y_\tau)d\tau$$
converges in probability to
$$\nu(t)\int_{-\ity}^{\ity}G^2(u)du$$
as $\epsilon \raro 0.$ This is a consequence of the central limit theorem for Ito stochastic integrals. (cf. Kutoyants (1984), p.78).
\vsp
\newsection{Proofs of Theorems}
\NI{\bf Proof of Theorem 2.1 :} From the inequality
$$(a+b+c)^2\leq 3(a^2+b^2+c^2), a,b,c\in R,$$
it follows that
\bea
\;\;\;\\\nonumber
E_\theta[|\hat \theta (t) X_t -\theta(t) x_t|^2]  &=& E_\theta [ |\frac{1}{\varphi_\epsilon}  \int^T_0 G \left(\frac{\tau-t}{\varphi_\epsilon} \right) \left(\theta(\tau) X_\tau -\theta(\tau) x_\tau \right)  d \tau \\ \nonumber
&& + \frac{1}{\varphi_\epsilon} \int^T_0 G \left(\frac{\tau-t}{\varphi_\epsilon}\right) \theta(\tau) x_\tau d \tau- \theta(t) x_t\\\nonumber
&& + \frac{\epsilon}{\varphi_\epsilon} \int^T_0 G \left(\frac{\tau-t}{\varphi_\epsilon} \right)
 \sigma_1(\tau,X_\tau,Y_\tau)dW_\tau|^2]\\ \nonumber
& \leq  & 3 E_\theta[ |\frac{1}{\varphi_\epsilon}  \int^T_0 G \left(\frac{\tau-t}{\varphi_\epsilon} \right) (\theta(\tau) X_\tau -\theta(\tau) x_\tau) d\tau|^2]\\ \nonumber
 & & + 3 E_\theta [|\frac{1}{\varphi_\epsilon} \int^T_0 G \left(\frac{\tau-t}{\varphi_\epsilon} \right)\theta(\tau) x_\tau d\tau -\theta(t)x_t |^2 ]\\ \nonumber
 & & +  3 \frac{\epsilon^2}{\varphi_\epsilon^2} E_\theta [ |\int^T_0 G \left( \frac{\tau-t}{\varphi_\epsilon}\right) \sigma_1(\tau,X_\tau,Y_\tau)dW_\tau|^2]\\ \nonumber
 &= & I_1+I_2+I_3. \;\;\\ \nonumber
\eea
By the boundedness condition on the function $\theta(.),$ the inequality (2.3) in Lemma 2.1  and applying the H\"older inequality, it follows that
\bea
\;\;\;\\\nonumber
I_1 &= &3 E_\theta \left| \frac{1}{\varphi_\epsilon} \int^T_0 G
\left(\frac{\tau-t}{\varphi_\epsilon} \right) (\theta(\tau) X_\tau -\theta(\tau) x_\tau)
d\tau \right|^2 \\\nonumber
&= & 3E_\theta  \left| \int^{\min(B,\frac{T-t}{\phi_\epsilon})}_{\max(A,-\frac{t}{\phi_\epsilon})} G(u) \left(\theta(t+\varphi_\epsilon u) X_{t+\varphi_\epsilon u}  - \theta(t+\varphi_\epsilon u) x_{t+\varphi_\epsilon u}\right) du\right|^2\\\nonumber
& \leq & 3 (B-A) \int^\infty_{-\infty} |G(u)|^2 L^2 E \left|X_{t+\varphi_\epsilon u}-x_{t+\varphi_\epsilon u} \right|^2 \ du
\;\;\mbox{(by using the condition $(A_1)$)}\\\nonumber
& \leq & 3(B-A)\int^\infty_{-\infty} |G(u)|^2 \;\;L^2 \sup_{0 \leq t +
\varphi_\epsilon u \leq T}E_\theta \left|X_{t+\varphi_\epsilon u}
-x_{t+\varphi_\epsilon u}\right|^2 \ du \\\nonumber
& \leq & 3 (B-A)L^2  e^{2LT} \epsilon^2 T \int_{-\ity}^\ity|G(u)|^2du\;\;\mbox{(by using the inequality (1.4))}\\\nonumber
\eea
and the last term tends to zero as $\epsilon \raro 0.$  By the boundedness condition on the function $\theta(.)$ and the H\"older inequality, it follows that
\bea
\;\;\;\\\nonumber
I_2 &= & 3E_\theta \left| \frac{1}{\varphi_\epsilon} \int^T_0 G\left(
\frac{\tau-t}{\varphi_\epsilon}\right) \theta(\tau) x_\tau d \tau - \theta(t) x_t\right|^2 \\ \nonumber
& = & 3  \left| \int^\infty_{-\infty} G(u)
\left(\theta(t+\varphi_\epsilon u) x_{t+\varphi_\epsilon u}-\theta(t) x_t \right)  \ du \right|^2
\\\nonumber
& \leq  & 3 (B-A)L^2 \varphi_\epsilon^2  \int_{-\ity}^\ity|uG(u)|^2 du\;\;\mbox{(by $(A_1)$)}.\\\nonumber
\eea
The last term  tends to zero as  $\varphi_\epsilon \rightarrow 0.$ We will now get an upper bound for the term $I_3.$ Note that
\bea
\;\;\;\\ \nonumber
I_3 &= & 3\frac{ \epsilon^2}{\varphi_\epsilon^2} E_\theta \left|\int^T_0 G \left(\frac{\tau-t}{\varphi_\epsilon}\right) \sigma_1(\tau,X_\tau,Y_\tau)dW_\tau\right|^2 \\ \nonumber
&=& 3 \frac{ \epsilon^2}{\varphi_\epsilon^2} \int_0^T (G \left(\frac{\tau-t}{\varphi_\epsilon}\right ))^2 E[\sigma(\tau,X_\tau,Y_\tau)]^2d\tau \\\nonumber
&\leq & C \frac{ \epsilon^2}{\varphi_\epsilon^2}\\\nonumber
\eea
for some positive constant $C$. This follows from the condition $(C1^\prime)$ and the Remarks made in the Appendix. Theorem 
2.1 is now proved by using the equations (3.1) to (3.4).\\

\NI {\bf Proof of Theorem 2.2 :} Let $J(t)=\theta(t) x_t.$ By the Taylor's formula, for any $x \in R,$
$$ J(y) = J(x) +\sum^k_{j=1} J^{(j)} (x) \frac{(y-x)^j}{j !} +[ J^{(k)} (z)-J^{(k)} (x)] \frac{(y-x)^k}{k!} $$
for some $z$ such that $|z-x|\leq |y-x|.$ Using this expansion, the equation (3.1) and the condition $(A_3)$  in the expression for $I_2$ defined in the proof of  Theorem 2.1, it follows that
\ben
\;\;\\\nonumber
I_2 & = & 3 \left[
\int^\infty_{-\infty} G(u) \left(J(t+\varphi_\epsilon u) - J(t) \right)  \ du \right]^2\\ \nonumber
&= & 3[ \sum^k_{j=1} J^{(j)} (t) (\int^\infty_{-\infty}G(u) u^j du )\varphi^j_\epsilon (j!)^{-1}\\\nonumber
& & \;\;\;\;+(\int^\infty_{-\infty}G(u) u^k (J^{(k)}(z_u) -J^{(k)} (x_t))du \;\varphi^k_\epsilon (k !)^{-1}]^2\\ \nonumber
\een
for some $z_u$ such that $|x_t-z_u|\leq |x_{t+\varphi_\epsilon u}-x_t| \leq C|\varphi_\epsilon u|.$ Hence
\bea
I_2 & \leq & 3  L^2 \left[  \int^\infty_{-\infty} |G(u)u^{k+1}|\varphi^{k+1}_\epsilon (k!) ^{-1}  du  \right]^2
\\ \nonumber
& \leq & 3 L^2 (B-A)(k!)^{-2} \varphi^{2(k+1)}_\epsilon\int^\infty_{-\infty} G^2(u) u^{2 (k+1)}\ du \\\nonumber 
&\leq & C \varphi_\epsilon^{2(k+1)}\\ \nonumber
\eea
for some positive constant $C$. Combining the equations (3.2), (3.4) and  (3.5), we get that there exists a positive constant $C$
such that 
$$ \sup_{a \leq t \leq b}E_\theta|\hat\theta_t X_t-\theta(t) x_t|^2 \leq C_3 (\epsilon^2 +  \varphi^{2(k+1)}_\epsilon +\epsilon^2 \varphi_\epsilon^{-2}). $$ 
Choosing $ \varphi_\epsilon = \epsilon^{\frac{1}{k+2}},$  we get that 
$$ \limsup_{\epsilon\rightarrow 0} \sup_{\theta(.) \in \Theta_{k+1} (L) } \sup_{a \leq t
\leq b} E_\theta|\theta(t)X_t - \theta(t)x_t|^2\epsilon^ {-\frac{2(k+1}{k+2}}< \infty. $$ This completes the proof of
Theorem 2.2. 
\vsp 
\NI{\bf Proof of Theorem 2.3:} Let $\alpha>0$ to be chosen later. From the equation (3.1), we obtain that
\bea
\;\;\;\;\\\nonumber
\lefteqn{\epsilon^{-\alpha}( \hat\theta(t)X_t -\theta(t)x_t)}\\\nonumber
 &= &\epsilon^{-\alpha}[\frac{1}{\varphi_\epsilon} \int^T_0 G \left(\frac{\tau-t}{\varphi_\epsilon} \right)
 \left( \theta(\tau)X_\tau-\theta(\tau) x_\tau\right) \  d \tau \\ \nonumber
 & & + \frac{1}{\varphi_\epsilon} \int^T_0 G \left( \frac{\tau-t}{\varphi_\epsilon}\right) \theta(\tau) x_\tau d\tau -\theta(t) x_t+ \frac{\epsilon}{\varphi_\epsilon} \int^T_0 G \left( \frac{\tau-t}{\varphi_\epsilon}\right)\sigma(\tau, X_\tau,Y_\tau)dW_\tau]\\\nonumber
& = &  \epsilon^{-\alpha}[\int_{-\ity}^{\ity} G(u)(\theta(t+\varphi_\epsilon u)X_{t+\varphi_\epsilon u}-\theta(t+\varphi_\epsilon u)x_{t+\varphi_\epsilon u})du]\\\nonumber
&&\;\;\;\;+\int_{-\ity}^{\ity}G(u)(\theta(t+\varphi_\epsilon u)x_{t+\varphi_\epsilon u}-\theta(t)x_t)du\\\nonumber
&&\;\;\;\;+\frac{\epsilon}{\varphi_\epsilon}\int_0^TG \left( \frac{\tau-t}{\varphi_\epsilon}\right)\sigma(\tau,X_\tau,Y_\tau)dW_tau\\\nonumber
 &=& R_1+R_2+R_3.\\\nonumber
\eea
By the boundedness  condition on the function $\theta(.)$ and the part (a) of Lemma 2.1, it follows that
\bea
R_1 & \leq & \epsilon^{-\alpha}|\int_{-\ity}^\ity G(u)(\theta(t+\varphi_\epsilon u) X_{t+\varphi_\epsilon u} - \theta(t+\varphi_\epsilon u)x_{t+\varphi_\epsilon u}) du|\\\nonumber
&\leq & \epsilon^{-\alpha} L \int_{-\ity}^\ity |G(u)| |X_{t+\varphi_\epsilon u} - x_{t+\varphi_\epsilon u}| du\\\nonumber
& \leq & Le^{LT} \epsilon^{1-\alpha} \int_{-\ity}^\ity |G(u)|\sup_{0\leq t+\varphi_\epsilon u \leq T}|\int_0^{t+\varphi_\epsilon u} \sigma_1(\tau,X_\tau,Y_\tau)dW_\tau|du\\\nonumber
\eea
by (3.3). Applying the Markov's inequality, it follows that, for any $\eta >0,$
\bea
\;\;\;\;\\\nonumber
P(|R_1|>\eta) &\leq &  \epsilon^{1-\alpha} \eta^{-1} Le^{LT} \int_{-\ity}^{\ity} |G(u)||E_\theta(\sup_{0\leq t+\varphi_\epsilon u \leq T} \int_0^{t+\varphi_\epsilon u}\sigma(\tau,X_\tau,Y_\tau)|du\\\nonumber
&\leq & \epsilon^{1-\alpha}\eta^{-1} Le^{LT} \int_{-\ity}^{\ity} |G(u)||E_\theta[(\sup_{0\leq t+\varphi_\epsilon u \leq T}\int_0^{t+\varphi_\epsilon u}\sigma(\tau,X_\tau,Y_\tau)dW_\tau)^2]|^{1/2} du\\\nonumber
&\leq & \epsilon^{1-\alpha}  \eta^{-1} Le^{LT} C T\int_{-\ity}^\ity|G(u)|du\\\nonumber
\eea
from the arguments given in Lemma 1.1, the condition $(A_1)$ and the properties of Ito stochastic integrals (cf. Prakasa Rao (1999), p. 31.).The last term tends to zero as $\epsilon \raro 0.$ Let $J_t=\theta(t)x_t.$ By the Taylor's formula, for any $t \in [0,T],$
$$ J_t = J_{t_0} + \sum^{k+1}_{j=1} J_{t_0}^{(j)}  \frac{(t-t_0)^j}{j !} + [ J_{t_0+\gamma(t-t_0)}^{(k+1)}-J_{t_0}^{(k+1)}] \frac{(t-t_0)^{k+1}}{(k+1)!} $$
where $0<\gamma<1$ and $t_0 \in (0,T).$ Applying the condition $(A_3)$ and the Taylor's expansion, it follows that
\bea
R_2 &=& \epsilon ^{-\alpha}[\sum_{j=1}^{k+1}J_t^{(j)}(\int_{-\ity}^\ity G(u) u^j \;du)\varphi_\epsilon^j(j!)^{-1}\\\nonumber
&& \;\;\;\; +\frac{\varphi_\epsilon^{k+1}}{(k+1)!}\int_{-\ity}^\ity G(u) u^{k+1}(J_{t+\gamma \varphi_\epsilon u}^{(k+1)}-J_t^{(k+1)})\;du]\\\nonumber
&=& \epsilon ^{-\alpha} \frac{J_t^{(k+1)}}{(k+1)!}\int_{-\ity}^\ity G(u) u^{k+1}\;du\\\nonumber
&& \;\;\;\; + \varphi_\epsilon^{k+1} \epsilon^{-\alpha}\frac{1}{(k+1)!}\int_{-\ity}^\ity G(u)u^{k+1}(J_{t+\gamma  \varphi_\epsilon u}^{(k+1)}-J_t^{(k+1)})\;du.\\\nonumber.
\eea
Observing that $\theta(t) \in \Theta_{k+1}(L),$ we obtain that
\bea
\lefteqn{\frac{1}{(k+1)!}\int_{-\ity}^\ity G(u)u^{k+1}(J_{t+\gamma\varphi_\epsilon u}^{(k+1)}-J_t^{(k+1)})du}\\\nonumber
&\leq & \frac{1}{(k+1)!}\int_{-\ity}^\ity |G(u)u^{k+1}(J_{t+\gamma\varphi_\epsilon u}^{(k+1)}-J_t^{(k+1)})|du\\\nonumber
&\leq & \frac{L\varphi_\epsilon}{(k+1)!}\int_{-\ity}^\ity |G(u)u^{k+2}|du.\\\nonumber
\eea
Combining the equations given above, it follows that
\bea
\lefteqn{\epsilon^{-\alpha} (\hat\theta_tX_t-J(t)- \frac{J_t^{(k+1)}}{(k+1)!} \int^\infty_{-\infty} G(u) u^{k+1}\ du)}\\\nonumber
&=& O_p(\epsilon^{1-\alpha})+O_p(\epsilon^{-\alpha}\varphi_\epsilon^{k+2})+\epsilon^{1-\alpha}\varphi_\epsilon^{-1}\int_0^TG(\frac{\tau-t}{\varphi_\epsilon}) \sigma(\tau,X_\tau,Y_\tau)dW_\tau.
\eea
Let
\be
\eta_\epsilon(t)= \epsilon^{-\alpha}\epsilon\varphi_\epsilon^{-1}\int_0^TG(\frac{\tau-t}{\varphi_\epsilon})\sigma(\tau,X_\tau,Y_\tau)dW_\tau.
\ee
Note that $E[\eta_\epsilon(t)]=0,$ by the properties of Ito integrals and 
\ben
E([\eta_\epsilon(t)]^2)&=& (\epsilon^{1-\alpha}\varphi_\epsilon^{-1})^2E([\int_0^TG(\frac{\tau-t}{\varphi_\epsilon})\sigma(\tau,X_\tau,Y_\tau)dW_\tau]^2)\\
&=& (\epsilon^{1-\alpha}\varphi_\epsilon^{-1})^2 \int_0^TG^2(\frac{\tau-t}{\varphi_\epsilon})E[\sigma(\tau,X_\tau,Y_\tau)]^2d\tau
\een
again by  the properties of Ito stochastic integrals as the processes $X$ and $Y$ are adapted to the filtration generated by the Brownian motion  $W.$ We choose $\alpha$ such that 
$$\epsilon^{1-\alpha} \varphi_\epsilon^{-1}=\varphi_\epsilon^{-1/2}.$$ 
Observe that $\alpha= \frac{2k+3}{2k+4}.$ Note that
\bea
\lefteqn{Var[\varphi_\epsilon^{-1/2}\int_0^TG \left(\frac{\tau-t}{\varphi_\epsilon} \right) \sigma(\tau,X_\tau,Y_\tau)dW_\tau]}\\\nonumber
&=&\varphi_\epsilon^{-1}\int_0^TG^2\left(\frac{\tau-t}{\varphi_\epsilon} \right) E[\sigma(\tau,X_\tau,Y_\tau)]^2d\tau\\\nonumber
&=&\int_{-\ity}^{\ity} G^2(u)\nu(t+\varphi_\epsilon u) du
\eea
and the last term tends to 
$$\nu(t)\int_{-\ity}^{\ity}G^2(u)du$$
as $\epsilon \raro 0$ where 
\be
\nu(\tau)=E[\sigma_1(\tau,X_\tau,Y_\tau)]^2.
\ee 
Observe that $\sup_{0\leq \tau \leq T}\nu(\tau) <\ity$ by the condition $(A_1)$ and the bounded convergence theorem. Following the equation (4.11) and the arguments given above, it follows that , as $\epsilon \raro 0,$  the asymptotic distribution of
$$ \epsilon^{-\frac{2k+3}{2k+4}} (\hat\theta_t X_t - J(t)  - \frac{J^{(k+1)}(t)}{(k+1) !} \int^\infty_{-\infty} G(u) u^{k+1}\ du)$$
is  the limiting distribution of the random variable
$$\frac{1}{\varphi_\epsilon}\int_0^T G \left(\frac{\tau-t}{\varphi_\epsilon} \right)\sigma_1(\tau,X_\tau,Y_\tau)dW_\tau$$
as $\epsilon \raro 0$ with mean zero and the variance
 $$\nu(t)\int_R G^2(u)du \;\;\mbox{where}\;\; \nu(t)=E[\sigma_1(t,X_t,Y_t)]^2.$$ 
This completes the proof of Theorem 2.3.
\vsp
\newsection{Alternate Estimator for the Multiplier $\theta(.)$}
Let $\Theta_\rho(L_\gamma)$ be a class of functions $\theta(.)$ uniformly bounded  by a constant $L$ and $k$-times continuously differentiable for some integer 
$k \geq 1 $ with the $k$-th derivative satisfying the H\"older condition of the order $\gamma \in (0,1):$
$$|\theta^{(k)}(t)-\theta^{(k)}(s)|\leq L_\gamma |t-s|^\gamma, \rho=k+\gamma.$$
Observe that $\rho >1.$ Suppose the process $\{X_t,0\leq t \leq T\}$ satisfies the stochastic differential equation given by the equation (1.1) where the linear multiplier is an unknown function in the class $\Theta_\rho(L_\gamma)$ and further suppose that $x_0 > 0$ and is {\it known}. 
From the Lemma 1.1, with probability one, it follows that
$$|X_t-x_t| \leq \epsilon e^{Lt}\sup_{0\leq v \leq T} |\int_0^v \sigma(s,X_s,Y_s)dW_s|.$$
Let
$$A_t= \{\omega: \inf_{0\leq s \leq t}X_s(\omega)\geq \frac{1}{2}x_0e^{-Lt}\}$$
and let $A=A_T.$ Following the technique suggested in Kutoyants (1994), p. 156, we define another process $Z$ with the differential
$$dZ_t=\theta(t) I(A_t) dt + \epsilon 2x_0^{-1} e^{LT} I(A_t)\sigma_1(t,X_t,Y_t)dW_t, 0\leq t \leq T.$$
We will now construct an alternate estimator of the linear multiplier $\theta(.)$ based on the process $Y$ over the interval $[0,T].$
Define the estimator
$$\tilde \theta(t)= I(A) \frac{1}{\varphi_\epsilon}\int_0^T G(\frac{t-s}{\varphi_\epsilon})dZ_s$$
where the kernel function $G(.)$ satisfies the condition $(A_2)$ and $(A_3)$. Observe that
\ben
E|\tilde \theta(t)-\theta(t)|^2 &= & E_\theta|I(A) \frac{1}{\varphi_\epsilon}\int_0^T G(\frac{t-s}{\varphi_\epsilon})(\theta(s)-\theta(t))ds\\\nonumber
&& \;\;\;\; + I(A^c)\theta(t)+I(A)\frac{\epsilon}{\varphi_\epsilon}\int_0^T G(\frac{t-s}{\varphi_\epsilon})2x_0^{-1}e^{LT}\sigma(s,X_s,Y_s)dW_s|^2\\\nonumber
&\leq & 3 E_\theta|I(A)\int_{-\ity}^{\ity} G(u)[\theta(t+u\varphi_\epsilon)-\theta(t)]du|^2+ 3 |\theta(t)|^2 [P(A^c)]^2\\\nonumber
&& \;\;\;\; + 3 \frac{\epsilon^2}{\varphi_\epsilon^2}|E[I(A)\int_0^T G(\frac{t-s}{\varphi_\epsilon})2x_0^{-1}e^{LT}\sigma(s,X_s,Y_s)dW_s]|^2\\\nonumber
&=& D_1+D_2+D_3. \;\;\mbox{(say)}.\\\nonumber
\een
Applying the Taylor's theorem and using the fact that the function $\theta(t)\in \Theta_\rho(L_\gamma)$, it follows that
\ben
D_1 \leq C_1\frac{1}{(k+1)!}\varphi_\epsilon^{2\rho} \int_{-\ity}^{\ity}|G^2(u)u^{2\rho}|du.
\een
Note that, by Lemma 1.1,
\ben
P(A^c) & = & P(\inf_{0\leq t \leq T}X_t < \frac{1}{2}x_0e^{-LT})\\\nonumber
&\leq & P(\inf_{0 \leq t \leq T}|X_t-x_t| + \inf_{0\leq t \leq T}x_t < \frac{1}{2}x_0e^{-LT})\\\nonumber
&\leq & P(\inf_{0 \leq t \leq T}|X_t-x_t| < -\frac{1}{2}x_0e^{-LT})\\\nonumber
&\leq & P(\sup_{0 \leq t \leq T}|X_t-x_t| > \frac{1}{2}x_0e^{-LT})\\\nonumber
&\leq & P(\epsilon e^{LT}\sup_{0 \leq t \leq T}|\int_0^t\sigma(s,X_s,Y_s)dW_s|>\frac{1}{2}x_0e^{-LT})\\\nonumber
&= & P(\sup_{0 \leq t \leq T}|\int_0^t\sigma(s,X_s,Y_s)dW_s|>\frac{x_0}{2\epsilon}e^{-2LT})\\\nonumber
&\leq & (\frac{x_0}{2\epsilon}e^{-2LT})^{-2}E[\sup_{0 \leq t \leq T}|\int_0^t\sigma(s,X_s,Y_s)dW_s|^2]\\\nonumber
&\leq & (\frac{x_0}{2\epsilon}e^{-2LT})^{-2}CT
\een
by the arguments given in Lemma 1.1 for some positive constant $C$ from the properties of Ito stochastic integrals (cf. Prakasa Rao (1999), p.31). The upper bound obtained above and the fact that $|\theta(s)|\leq L, 0\leq s \leq T$ leads  an upper bound for the term $D_2.$ We have used the inequality
$$x_t= x_0 \exp(\int_0^t\theta (s)ds)\geq x_0 e^{-Lt}$$
in the computations given above. Applying Theorem 2.1, it follows that
\ben
\lefteqn{E[|I(A)\int_0^T G(\frac{t-s}{\varphi_\epsilon})2x_0^{-1}e^{LT}\sigma(s,X_s,Y_s)dW_s|^2]}\\\nonumber
&\leq & CE[|\int_0^T G(\frac{t-s}{\varphi_\epsilon})\sigma(s,X_s,Y_s)dW_s|^2]\\\nonumber
&=& C\;Var[\int_0^T G(\frac{t-s}{\varphi_\epsilon})\sigma(s,X_s,Y_s)dW_s]\\\nonumber
&=& C\varphi_\epsilon \nu(t)\int_RG^2(u) du +o(1)\\\nonumber
\een
as $\epsilon \raro 0$ for some positive constant $C$ which leads to an upper bound on the term $D_3.$ Combining the above estimates, it follows that
\ben
E|\tilde \theta(t)-\theta(t)|^2\leq C(\varphi_\epsilon^{2\rho} + \epsilon^4+ \epsilon^2 \varphi_\epsilon)
\een
for some positive constants $C.$ Choosing $\varphi_\epsilon=\epsilon^{\frac{2}{2\rho-1}},$ we obtain that
\ben
E|\tilde \theta(t)-\theta(t)|^2\leq C (\epsilon^{\frac{4\rho}{2\rho-1}}+  \epsilon^{4})\leq C \epsilon^{\frac{4\rho}{2\rho-1}}
\een
for some positive constant $C$ since $\rho >1.$ Hence we obtain the following result implying the uniform consistency of the estimator $\tilde \theta(t)$ as an estimator of $\theta(t)$ as $\epsilon \raro 0.$ 
\vsp
\NI{\bf Theorem 4.1:} {\it Let $\theta \in \Theta_\rho(L)$. Let $\varphi_\epsilon= \epsilon^{2/(2\rho-1)}.$ Suppose the conditions $(A_1)-(A_3), (C1^\prime)$ and $(C2^\prime)$ hold. Then, for any interval $[a,b] \subset [0,T],$ }
\ben
\limsup_{\epsilon \raro 0}\sup_{\theta(.)\in \Theta_\rho(L)}\sup_{a\leq t \leq b}E|\tilde \theta(t)-\theta(t)|^2 \epsilon^{- 
\frac{4\rho}{2\rho-1}}<\ity.
\een
\vsp
\CC{\bf Appendix} 
Khlifa et al. (2017) studied a class of stochastic differential equations (SDEs) of the form 
$$
dX_t=a(t,X_t)dt+\sigma(t,X_t,Y_t)\;dW_t, X_0=x_o\in R, 0\leq t \leq T
$$
where $\{W_t, t\geq 0\}$ is a standard Brownian motion and $\{Y_t, t \geq 0\}$ is a continuous stochastic process adapted to the filtration generated by the Brownian motion $W.$ They proved the following result concerning the existence and the uniqueness of the solution for such SDEs.
\vsp
\NI{\bf Theorem :} {\it Let $a(.,.)$ and $\sigma(.,.,.)$ be non-random functions and $Y$ be an adapted continuous process. Suppose the following conditions hold:\\
(C1) there exists a constant $K>0$ such that
$$
|\sigma(t,x,y)|^2+|a(t,x)|^2\leq K^2(1+|x|^2), t \geq 0, x,y \in R, 
$$
(C2) for any integer $N\geq 1,$ there exists a constant $ K_N>0$ such that, for all $t \geq 0$ and for all $x_1,x_2,y$ satisfying $|x_1|\leq N, |x_2|\leq N$ and $|y|\leq N,$
$$
|a(t,x_1)-a(t,x_2)|\leq K_N|x_1-x_2|
$$
and
$$
|\sigma(t,x_1,y)-\sigma(t,x_2,y)|\leq K_N \phi(t,y)|x_1-x_2|
$$
where $\phi(.,.)$ is a positive function such that for any integer $N \geq 1,$ there exists a constant $ C_N>0$ such that 
$$
\sup_{s \geq 0}\sup_{|x|\leq N}|\phi(s,x)|\leq C_N.
$$
Then the equation 
$$dX_t=a(t,X_t)dt+\sigma(t,X_t,Y_t)\;dW_t, X_0=x_o\in R, 0\leq t \leq T$$
has a unique strong solution.}
\vsp
\NI{\bf Remarks:} Suppose the conditions stated in the theorem given above hold. Following Lemma 1 of Khlifa et al.(2016) and the arguments in Skorokhod (1965), it can be shown that the process $\{X_t,0\leq t \leq T\}$ has the property (cf. Theorem 5.4, Klebaner (2012)) that
$$
E(\sup_{0\leq t \leq T}X_t^2)\leq C( 1+E(X_0^2))
$$ 
for some constant $C>0.$ 
\vsp
{\bf Funding support:} There is no special funding for this work except for the fact that this work was done under the scheme ``INSA Honorary Scientist" supported by the Indian National Science Academy at the CR RAO Advanced Institute of Mathematics, Statistics and Computer Science, Hyderabad 500046, India .
\vsp
{\bf Conflict of Interest:} There is no conflict of interest and hundred percent of the work reported here is done by me.
\vsp
\NI{\bf References :}
\vsp
\begin{description}
\item Ait-Sahalia and Kimmel, R. (2007) Maximum likelihood estimation of stochastic volatitlity models, {\it J. Financ. Econ.}, {\bf 83}, 413-452.
\item Fouque, J.-P., Papanicolaou, G., and Sircar, K.R. (2000a) {\it Derivatives in Financial Markets with Stochastic Volatility}, Cambridge University Press, Cambridge.
\item Fouque, J.-P., Papanicolaou, G., and Sircar, K.R. (2000b) Mean-reverting stochastic volatility, {\it Int. J. Theor. appl. Finance}, {\bf 3}, 101-142.
\item Heston, S. (1993) A closed-form solution of options with stochastic volatility with applications to bond and currency options, {\it The Review of Financial Studies}, {\bf 6}, 327-343.
\item Hull, J. and White, A. (1987) The pricing of options on assets with stochastic volatilities, {\it J. Finance.}, {\bf 42}, 281-300.
\item Khlifa, M.B.H., Mishura,Y., Ralchenko, K., and Zili, M. (2016) Drift parameter estimation in stochastic differential equation with multiplicative stochastic volatility, {\it Modern Stochastics: Theory and Methods}, {\bf 3}, 269-285.
\item Khlifa, M.B.H., Mishura,Y., Ralchenko, K., Shevchenko, G., and Zili, M. (2018) Stochastic differential equations with generalized stochastic volatility and statistical estimators, {\it Theor. Probability and Math. Statist.}, {\bf 96}, 1-13.
\item Klebaner, F.C. (2012) {\it Introdcution to Stochatic Calculus with Applications}, Imperial College Press, London.
\item Kuchuk-Iatsenko,S. and Mishura, Y. (2015a) Option pricing in the model with stochastic volatility driven by Ornstein-Uhlenbeck process. Simulation., {\it Mod. Stoch. Theory Appl.}, {\bf 2}, 355-369.
\item Kuchuk-Iatsenko,S. and Mishura, Y. (2015b) Pricing the European call option in the model with stochastic volatility driven by Ornstein-Uhlenbeck process, {\it Mod. Stoch. Theory Appl.}, {\bf 3}, 233-249.
\item Kutoyants, Y. (1984) {\it Parameter Estimation for Stochastic Processes}, Translated and Edited by B.L.S. Prakasa Rao, Helderman Verlag, Berlin.
\item Kutoyants, Y. (1994) {\it Identification of Dynamical Systems with Small Noise}, Kluwer Academic Publishers, Dordrecht.
\item Prakasa Rao, B.L.S. (1999) {\it Semimartingales and Their Statistical Inference}, Chapman and Hall and CRC Press, Boca Raton.
\item Skorokhod, A.V. (1965) {\it Studies in the Theory of Random Processes}, Addison-Wesley Publishers, Reading, Mass.

\end{description}

\end{document}